\documentclass{gen-j-l}

\newtheorem{theorem}{Theorem}

\newtheorem{statement}[theorem]{Statement}

\theoremstyle{definition}

\newtheorem{example}[theorem]{Example}

\theoremstyle{remark}

\numberwithin{equation}{section}

\begin{document}
  \title{A correction to a result of B.~Maier}
  \author[A. Magidin]{Arturo Magidin}
  \address{Oficina 112\\
           Instituto de Matem\'aticas\\
           Universidad Nacional Aut\'onoma de M\'exico\\
           Circuito Exterior\\
           Ciudad Universitaria\\
           04510 Mexico City, M\'exico\\}
  \email{magidin@matem.unam.mx}


\begin{abstract}
In a 1985 paper, Berthold J.~Maier gave necessary and sufficient
conditions for the weak embeddability of amalgams of two nilpotent groups
of class two over a common subgroup. Then he derived simpler
conditions for some special cases. One of his subsequent results is
incorrect, and we provide a counterexample. Finally, we provide a fix
for the result.
\end{abstract}

\maketitle

In this note, $\mathcal{N}_2$ denotes the variety of nilpotent groups
of class two, that is, groups $G$ such that $[G,G]\subseteq
Z(G)$. Recall that an \textit{amalgam of $A$ and $B$ over the common
subgroup $D$} consists of two groups, $A$ and $B$, and a group $D$
which is a subgroup of both $A$ and $B$. The amalgam is \textit{weakly
embeddable} in $\mathcal{N}_2$ if and only if there is a
$\mathcal{N}_2$-group $G$, such that $A$ and $B$ are subgroups of $G$
and $D\subseteq A\cap B$ (inside~$G$). We then say that $G$ is a
(weak) amalgam for $A$ and $B$ over $D$. If $G$ satisfies the further
condition that $D= A\cap B$, then $G$ is said to be a \textit{strong}
amalgam.

In \cite{amalgone}, Berthold J.~Maier studied the question of weak
embeddability of amalgams in $\mathcal{N}_2$. Note that when he says
that an amalgam exists, he means that a \textit{weak} amalgam
exists. The Hauptsatz in \cite{amalgone} is a characterization of weak
embeddability for amalgams in $\mathcal{N}_2$. We quote it here for reference:

\begin{theorem}[Berthold J.~Maier, Haupsatz in \cite{amalgone}]
Let $A,B\in\mathcal{N}_2$, with a common subgroup $D\leq A,B$. There
exists a weak amalgam of $A$ and $B$ over $D$ in $\mathcal{N}_2$ if and only
if the following two conditions hold:
\begin{itemize}
\item [(1)] $A_2\cap D\leq Z(B)$ and $B_2\cap D\leq Z(A)$.
\item[(2)] For all $k>0$, $q_i>0$, $x_i\in A$ and $x'_i\in A_2$ with
$x_i^{q_i}x'_i\in D$, $y_i\in B$ and $y'_i\in B_2$ with
$y_i^{q_i}y'_i\in D$, we have that for every $d\in D$,
\[\prod_{i=1}^k \left[x_i,y_i^{q_i}y'_i\right] = d \Longleftrightarrow 
  \prod_{i=1}^k \left[x_i^{q_i}x'_i,y_i\right] = d.\]
\end{itemize}
\end{theorem}

In the above, $A_2$ is the commutator subgroup of $A$, $Z(A)$ is the center of
$A$, and similarly for $B$; the commutator brackets are given by
$[x,y]=x^{-1}y^{-1}xy$. 

After the proof, Maier investigates whether the conditions may be
relaxed given more information about $A$ and $B$. For example, he
proves that if $A$ and $B$ are torsion-free, then condition (2) above
is superfluous. He also shows that if $D$ is normal in either $A$ or
$B$, then it suffices to consider condition (2) with $k=1$. 

Then Maier looks at the case where $D$ is either central or
\textit{co-central} in~$B$. Maier defines a subgroup $X$ of $B$ to be
co-central if there exists a central subgroup $Z$ such that $\langle
X,Z\rangle = B$. The case where $D$ is central or co-central in $A$ is
covered by symmetry.

Maier's result reads:

\begin{statement}[Satz 2 in \cite{amalgone}]
\label{wrongresult} Let $A$ and $B$ in
$\mathcal{N}_2$ with a common subgroup $D$.

\begin{enumerate}
\item If $D$ is co-central in $B$, then there exists an amalgam of $A$
and $B$ over $D$ in $\mathcal{N}_2$ if and only if
\begin{itemize}
\item[(*)] For $q>0$, $a\in A$ with $a^q\in A_2D$, and $b\in B\setminus D$ with
$b^q\in D$, we have $[a,b^q]=e$.
\end{itemize}

\item If $D$ is central in $B$, then there exists an amalgam of
$A$ and $B$ over $D$ in $\mathcal{N}_2$ if and only if $B_2\cap
D\subseteq Z(A)$, and for $q>0$, $a\in A$ with $a^q\in A_2D$, and
$b\in B$, $b'\in B_2$ with $b^qb'$, we have $[a,b^qb']=e$.
\end{enumerate}
\end{statement}

In fact, Part 1 of the statement above is incorrect. Although the
condition given is sufficient, as Maier proves, it is not
necessary. Part 2, on the other hand, is correct in its entirety.

First, we provide a counterexample to the statement in Part 1:

\begin{example}
Two groups, $A$ and $B$ in $\mathcal{N}_2$, with a common subgroup $D$,
$D$ co-central in $B$, which fails (*), and such that the amalgam of
$A$ and $B$ over $D$ is embeddable in $\mathcal{N}_2$.

Let $D$ be the relatively free $\mathcal{N}_2$-group in two generators,
$x$ and $y$. The elements of the group can be uniquely written in the
form
\[x^a y^b [y,x]^c;\qquad a,b,c\in\mathbb{Z}.\]

Let $Z$ denote the infinite cyclic group, and $Z/qZ$ the cyclic group
of $q$ elements. We denote the generator of the group by $z$.

Let $A=B=Z/qZ\times D$. Clearly, $Z/qZ \times D \times Z/qZ$ is a
(strong) amalgam of $A$ and $B$ over $D$. Also, since $B=\langle D,
Z/qZ\rangle$, and the $Z/qZ$ factor of $B$ is central, $D$ is
co-central in $B$.

Let $a=(e,x)\in A$, and let $b=(z,y)\in B\setminus D$. Then
$a^q=(e,x^q)\in D$, and $b^q=(e,y^q)\in D$. However, 
\[ [a,b^q] = [x,y^q] = [x,y]^q\not=e. \] 
\end{example}

\bigskip

An easy way to fix the problem is to modify condition (*). Note first
that a subgroup $H$ of $G$ is co-central if and only if $\langle H,
Z(G) \rangle$ equals the whole group. To fix Maier's result, instead
of asking that $b\in B\setminus D$, we request merely that it lie in
the center of $B$.

To prove that this modification will work, we need to quote another
result of Maier's:

\begin{theorem}[Korollar 3 in \cite{amalgone}]
\label{fornormal}
Let $A$ and $B$ in $\mathcal{N}_2$ with a common subgroup $D$, and assume
that $D$ is normal in either $A$ or $B$. Then there exists an amalgam
of $A$ and $B$ over $D$ in $\mathcal{N}_2$ if and only if the following
two conditions hold:
\begin{itemize}
\item[(a)] $A_2\cap D\subseteq Z(B)$ and $B_2\cap D\subseteq Z(A)$.
\item[(b)] For $q>0$, $a\in A$, $a'\in A_2$ with $a^qa'\in D$, $b\in
B$, $b'\in B_2$ with $b^qb'\in D$, we have 
\[[a^qa',b] = [a,b^qb'] \in D.\]
\end{itemize}
\end{theorem}

We can now provide a corrected version of Statement~\ref{wrongresult},
Part 1:

\begin{theorem}
Let $A$ and $B$ in $\mathcal{N}_2$, and let $D$ be a common subgroup of
$A$ and $B$. Assume further that $D$ is co-central in $B$.
There exists a (weak) amalgam of $A$ and $B$ over $D$
in $\mathcal{N}_2$ if and only if the following condition holds:
\begin{itemize}
\item[(**)] For $q>0$, $a\in A$ with $a^q\in A_2D$, and $z\in Z(B)$ with
$z^q\in D$, we have $[a,z^q]=e$.
\end{itemize}
\end{theorem}

\begin{proof}
First we prove necessity. Let $G$ be an amalgam of $A$ and $B$ over
$D$, with $D$ co-central in $B$. Let $q$, $a$, and $z$ be as in
(**). Since $a^qa'\in D\subseteq B$, and $z\in Z(B)$, $[a^qa',z]=e$. But
in $G$ the commutator bracket acts bilinearly on $G^{{\rm ab}}\times
G^{{\rm ab}}$, so
\begin{eqnarray*}
[a,z^q] & = & [a^q,z]\\
& = & [a^qa',z]\\
& = & e.
\end{eqnarray*}
This proves (**), giving necessity.

We will prove sufficiency by showing that the hypothesis that $D$ is
co-central in $B$ together with (**) imply (a) and (b) from
Theorem~\ref{fornormal}. Note that a co-central subgroup is
necessarily normal, so Theorem~\ref{fornormal} applies.

To prove (a), let $a\in A_2\cap D$. Then $a$ is necessarily in $Z(D)$,
since $A$ lies in $\mathcal{N}_2$; thus $a$ commutes with every element
of $D$, and when considered inside of $B$, it must commute with
everything in the center of $B$. Since $B$ is generated by $D$ and
$Z(B)$, it follows that $a\in Z(B)$. Thus, $A_2\cap D\subseteq Z(B)$.

For the other inclusion, note that since $D$ is co-central and $B$ is
nilpotent of class two, it follows that $B_2=D_2$. In particular,
\[B_2\cap D = D_2 \subseteq A_2\subseteq Z(A).\]

To prove (b), let $a\in A$, $a'\in A_2$, $b\in B$, $b'\in B_2$, and
$q>0$ such that $a^qa'$ and $b^qb'$ both lie in $D$. We want to show that
\[[a^qa',b] = [a,b^qb'] \in D.\]

Since $B_2=D_2$, note that $b^qb'\in D$ if and only if $b^q\in D$, and
that $[a,b^qb']=[a,b^q]$ in $A$. So we may assume without loss of
generality that $b'=e$. 

Write $b=dz$, with $d\in D$ and $z\in Z(B)$. Then $b^q=d^qz^q$, so
$b^q\in D$ implies that $z^q\in D$. Then condition (**) applies to $a$
and $z$, so we conclude that $[a,z^q]=e$.

Therefore, in $A$, we have:
\begin{eqnarray*}
[a,b^q] & = & [a,(dz)^q]\\
        & = & [a,d^qz^q]\\
        & = & [a,d^q][a,z^q]\\
        & = & [a,d^q]\\
        & = & [a^q,d]\\
        & = & [a^qa',d].
\end{eqnarray*}
The third step can be done because both $d^q$ and $z^q$ lie in
$D$. The next to last and last steps follow because $A\in\mathcal{N}_2$.

On the other hand, in $B$ we have:
\begin{eqnarray*}
[a^qa',b] & = & [a^qa',dz]\\
          & = & [a^qa',d][a^qa',z]\\
          & = & [a^qa',d]
\end{eqnarray*}
since $z$ is central in $B$. 

Therefore, $[a,b^q] = [a^qa',b]$. Moreover, since
$[a^qa',b]\in B_2=D_2\subseteq D$, the two commutators lie in $D$. This
proves condition (b), and thus the theorem.
\end{proof}

\bibliographystyle{amsplain}
\bibliography{bibliog}

\providecommand{\bysame}{\leavevmode\hbox to3em{\hrulefill}\thinspace}
\begin{thebibliography}{1}

\bibitem{amalgone}
Berthold~J. Maier, \emph{Amalgame nilpotenter {G}ruppen der {K}lasse zwei},
  Publ. Math. Debrecen \textbf{31} (1985), 57--70, {MR {\bf 85k}:20117}.

\end{thebibliography}

\end{document}